\documentclass[a4paper,10pt]{amsart}
\usepackage{graphicx}
\usepackage{amssymb}
\usepackage[cp1251]{inputenc}
\usepackage[T2A]{fontenc}

\newtheorem{theorem}{Theorem}

\newtheorem{definition}{Definition}
\newtheorem{example}[definition]{Example}
\newtheorem{proposition}[definition]{Proposition}
\newtheorem{corollary}[definition]{Corollary}

\newcommand{\Q}{\mathbb{Q}}
\newcommand{\Z}{\mathbb{Z}}
\newcommand{\R}{\mathbb{R}}
\newcommand{\Co}{\mathbb{C}}

\newcommand{\Cst}{C^{*}}
\newcommand{\Img}{\emph{Im}}
\newcommand{\lprod}[3]{\;_{#1}\left< #2 , #3 \right>}
\newcommand{\rprod}[3]{\left< #2 , #3 \right>_{#1}}
\newcommand{\iprod}[2]{\left< #1 , #2 \right>}
\newcommand{\norm}[1]{||#1||}
\newcommand{\mat}[1]{\mathcal{M}#1}
\newcommand{\Mat}[2]{\mathcal{M}_{#1}#2}
\newcommand{\Siegel}[1]{\mathbb{H}_{#1}}

\title[The graded ring of quantum theta functions]{The graded ring of quantum
theta functions for noncommutative torus with real multiplication}

\author{Mariya Vlasenko}

\email{mariyka@imath.kiev.ua}

\begin{document}
\maketitle

\begin{abstract}
For quantum torus generated by unitaries $UV = e(\theta)VU$ there
exist nontrivial strong Morita autoequivalences in case when $\theta$ is real
quadratic irrationality. A.Polishchuk introduced and studied the graded ring of
holomorphic sections of powers of the respective bimodule (depending on the
choice of a complex structure). We consider a Segre square of this ring
whose graded components are spanned by Rieffel scalar products of Polishchuk's
holomorphic vectors as in~\cite{Boca1999} and~\cite{Manin2003}. These graded
components are linear spaces of quantum theta functions in sense of
Yu.Manin.
\end{abstract}

\section*{Introduction}
A quantum torus $A_{\theta}$ with an irrational parameter $\theta \in
\R\backslash\Q$ is a transformation group $\Cst$-algebra
$\Cst(\theta\Z,\R\slash\Z)$ for the group action of $\theta \Z$ on
$\R\slash\Z $ or, equivalently, a universal $\Cst$-algebra generated by two unitaries $U,V \in
A_{\theta}$ satisfying relation $UV = e(\theta)VU$. Here $e(x) = exp(2
\pi i x)$.

\begin{definition} $A_{\theta}$ is a quantum torus with real
multiplication if
$\theta$ is a real quadratic irrationality, i.e. a real irrational root of a
quadratic equation with rational coefficients.
\end{definition}

Let $k$ be a real quadratic field. In \cite{AG0202109} it is proposed to use
quantum tori with real multiplication $A_{\theta}$, $\theta \in k \backslash \Q$
as geometric objects associated to $k$. This should be compared to
consideration of elliptic curves with complex multiplication $E_{\tau} = \Co
\slash \Gamma$, $\Gamma = \Z + \tau \Z$, $\tau \in k^{'} \backslash \Q$
for complex quadratic field $k^{'}$. Any endomorphism $\alpha : E_{\tau}
\rightarrow E_{\tau}$ is a linear map on the universal covering $\Co$, so
$End(E_{\tau})$ is identified with the ring of multipliers of the lattice
$\Gamma$, that is $\{ \alpha \in \Co | \alpha \Gamma \subset \Gamma\}$. We say
that $E_{\tau}$ is an elliptic curve with complex multiplication
if $End(E_{\tau})$ is larger then $\Z$, which happens precisely when $\tau$ is
a complex quadratic number. 

Real multiplication of quantum tori has similar interpretation when we consider
morphisms in sense of noncommutative geometry: every element of $End
(A_{\theta})$ is by definition an (isomorphism class of)
$A_{\theta}$-$A_{\theta}$-bimodule, finitely generated and projective as
left and right module at the same time. Every such isomorphism class $[M] \in
End(A_{\theta})$ defines an endomorphism $\phi_{[M]}$ of $K_0$-group of
$A_{\theta}$ via $[P] \mapsto [P \underset{A_{\theta}}\otimes M]$ for finitely
generated projective
right $A_{\theta}$-modules $P$. It is shown in \cite{AG0202109} that when
$K_0(A_{\theta})$ is identified with the lattice $\Gamma = \Z + \theta \Z$
via trace map, then $\phi_{[M]}$ becomes a multiplication by real number.
Moreover, this map
$$
K_0 : End(A_{\theta}) \rightarrow \{ \alpha \in \R | \alpha \Gamma \subset
\Gamma\},\; K_0([M]) = \phi_{[M]}
$$ 
is surjective. So, $A_{\theta}$ is a quantum tori with real multiplication if
and only if $K_0(End(A_{\theta}))$ is larger then $\Z$.

In this paper we construct the graded ring of quantum theta functions $R
=\oplus R_n$ for quantum torus with real multiplication $A_{\theta}$. The
construction is described in Section~6, were we also prove that $\theta \in k
\backslash \Q$ can be chosen such that the ring $R$ is generated over $\Co$ by
finite dimensional vector space $R_1$ (Theorem~\ref{fin_gen}). 

We sketch the definition of $R$ below. First, we need simple facts from
number theory. One can prove that $\{ \alpha \in \R | \alpha \Gamma \subset
\Gamma\} = \Z + f O_k$ for some integer $f \ge 1$, where $O_k$ is the ring of
integers of real quadratic field $k = \Q(\theta)$. Thus there are units of infinite 
order of $O_k$ in $\Z + f O_k$, and we take one of them $\varepsilon \in (\Z + f O_k)
\cap O^{\times}_k$. Then there exists a bimodule $M_{\varepsilon}$ with
$K_0([M_{\varepsilon}])=\varepsilon$, which is an
$A_{\theta}-A_{\theta}$-imprimitivity bimodule. Bimodules of such a kind were
studied in \cite{QA0203160},\cite{QA0211262},\cite{AG0202109}, and we describe
them in Section~5. Sections~1-3 contain discussion of relationship between
biprojective bimodules and imprimitivity bimodules. $M_{\varepsilon}$ is
an infinite dimensional $\Co$-vector space, but one can take finite dimensional subspaces
$E_n \subset M_{\varepsilon}^{\otimes n}$ of so-called ``holomorphic'' vectors,
and they are compatible with tensor product: $E_n \otimes E_m \subset E_{n+m}$
(\cite{QA0203160},\cite{QA0211262}). So, we obtain the graded ring $E =
\oplus_n E_n$ with multiplication defined via tensor product.
This ring was studied in \cite{AG0212306}. The choice of ``holomorphic'' vectors
depend on a complex parameter $\tau \in \Co$, formally defining 
``holomorphic'' structure on $A_{\theta}$. In this paper we use structure of
imprimitivity bimodules on
$M_{\varepsilon}^{\otimes n}$ to obtain quantum theta functions from
``holomorphic'' vectors (see Section~4 for definition of the structure of
imprimitivity bimodule on tensor product). It was already noticed in
\cite{Boca1999} that operator-valued theta functions appear from imprimitivity
bimodules over quantum tori.

We use the definition of quantum theta functions given in~\cite{Manin1990} and~\cite{AG0011197}. Let us briefly recall it. Consider Heisenberg group $G_{\theta}$
$$
1 \rightarrow \Co^{\times} \rightarrow G_{\theta} \rightarrow \Co^2 \times \Z^2
\rightarrow 0
$$
acting on elements of quantum torus $A_{\theta}$ by
$$
(\alpha; \vec{x};\vec{m})\sum_{\vec n \in \Z^2} a_{\vec n} U^{n_1}V^{n_2} =
 \alpha \sum_{\vec n \in \Z^2} e(n_1 x_1 + n_2 x_2)a_{\vec n}
U^{m_1}V^{m_2}U^{n_1}V^{n_2}.
$$
A multiplier $\mathcal{L}$ is any free subgroup of rank 2 in $G_{\theta}$, which
is a lift of a free subgroup of rank 2 in $\Co^2 \times \Z^2$. We denote by
$\Gamma(\mathcal{L}) \subset A_{\theta}$ the vector space of elements fixed by
$\mathcal{L}$. All elements of $\Gamma(\mathcal{L})$ are called quantum
theta functions with multiplier $\mathcal{L}$. For example, take a lattice
$L = \Z \vec s + \Z \vec r \subset \Z^2$, and a matrix $\Omega \in \Mat 2
{\Co}$, symmetric $\Omega = \Omega^t$ and with positive imaginary part $\Im
\Omega > 0$. Then $(e(\frac 12 \vec s^t
A^t \vec s);A \vec s; \vec s)$ and $(e(\frac 12 \vec r^t A^t \vec r);A \vec r;
\vec r)$ generate a multiplier, where $A = \frac {\theta} 2
\begin{pmatrix}0&1\\-1&0\\\end{pmatrix} + \Omega$. Let us denote this multiplier
by $\mathcal{L} = \mathcal{L}(L, \Omega)$. Then $\Gamma(\mathcal{L})$ is
$\#(\Z^2 \slash L)$-dimensional $\Co$-vector space of elements of the form
$$
\Theta[f](\Omega) = \sum_{\vec{m} \in \Z^2} f(\vec m)e(\frac 12 \vec m ^t
\Omega \vec m)e(-\frac{\theta}2 m_1 m_2)U^{m_1}V^{m_2} 
$$ 
for $f: \Z^2 \slash L \rightarrow \Co$. 

We define graded components of our ring of quantum theta functions
as $R_n = \Gamma(\mathcal{L}(c_n \Z^2, \Omega_n)) \subset A_{\theta}$. Now
we explain what are $c_n$ and $\Omega_n$. Recall we have chosen a unit
$\varepsilon \in O_k$ and a complex parameter $\tau \in \Co$. Now we need them
to satisfy some technical conditions, especially $\Im \tau > 0$ and $\varepsilon
= c \theta + d > 0$, $\varepsilon \theta = a
\theta + b$ with $\begin{pmatrix}a&b\\c&d\\\end{pmatrix} \in SL_2(\Z)$ and $c
> 0$. Existence of such an $\varepsilon$ one can get for example from
Section~6. Then $c_n$ are defined by $\varepsilon^n = c_n \theta + d_n$
with integer $c_n, d_n$, or, equivalently, by 
$\begin{pmatrix}a&b\\c&d\\ \end{pmatrix}^n =
\begin{pmatrix}a_n&b_n\\c_n&d_n\\ \end{pmatrix}$, or by
$\sum_n c_n t^n = \frac c {t^2 -(a+d)t + 1}$. In particular, the last expression
shows that $\{ c_n \}$ is an increasing sequence of positive integers,
since $a+d = \varepsilon + \frac 1{\varepsilon} \ge 2$.
Now $\Omega_n = \frac {1}{c_n
\varepsilon^n} \Omega$ where  
$$
\Omega = \frac{i}{2 \Im
\tau}\begin{pmatrix}|\tau|^2&-\Re \tau\\-\Re \tau & 1 \\\end{pmatrix}.
$$

If we whish to define a product on $R = \oplus R_n$, then usual product in
$A_\theta$ wouldn't do. It was already noticed in \cite{AG0011197} that a
product of two quantum theta functions in $A_{\theta}$ is not a quantum
theta function as a rule. But our quantum theta functions in $R_n$
arise from bimodules $M^n_{\varepsilon}$
(Proposition~\ref{quantum_theta}), so we get another product~--- a bilinear
operation $\star: R_n \otimes R_m \rightarrow R_{n+m}$, naturally coming from
tensor product of bimodules.

So, we get the ring $R$ whose elements are formally an elements of
quantum tori $A_{\theta}$, but multiplication law is different from the one in
$A_{\theta}$. In fact $R$ is isomorphic to a kind of Segre square of $E$~--- the
subspace in $E \bar \otimes E$ generated by elements $a \bar \otimes b$ with
$a,b \in E_n$ for some $n$. Here $\bar \otimes$ means that
$(\alpha a)\bar\otimes b = a \bar\otimes(\bar \alpha b)$ for $\alpha \in \Co$.
Both $R$ and $E$ encapsulate the structure of real multiplication and use
arithmetical data to be constructed. But the following question still remains
unanswered: whether we can use such a rings to obtain arithmetical invariants of
real quadratic field  $k = \Q(\theta)$?  
\\
\\
Acknowledgement. I am grateful to Yuri Ivanovich Manin who taught me the idea
of real multiplication and collaboration with whom led to writing this
notes. I also thank Alexander Polishchuk who read 
the first version of the paper and remarked that our ring of quantum theta
functions is a Segre square of the ring of ``holomorphic'' vectors.

\numberwithin{definition}{section}
\section{Strong Morita equivalence}

Let $A$ be a pre-$\Cst$-algebra, i.e. $\Co$-algebra with involution and
norm satisfying $\norm{x}^2=\norm{x^*x}$ and $\norm x = 0$ if and only if $x =
0$ for $x \in A$. If $A$ has a unit element $1 \in A$ it is assumed that $\norm
1
= 1$. An $A$-valued pre-inner product on linear space $M$ is an $A$-valued
sesquilinear form $\iprod{\cdot}{\cdot}$ (here it does not matter in which
variable it is conjugate linear) such that $\iprod x x \ge 0$ in completion of
$A$ and $\iprod x y ^*= \iprod y x$ for $x,y \in M$. We denote
$\Img \iprod{\cdot}{\cdot} \subset A$ the set of finite sums of elements of
the form $\iprod x y$ for $x,y \in M$.

Following definitions were introduced in \cite{Riff1974}.

\begin{definition}
A left $A$-module $M$ is called a left $A$-rigged space if it is
endowed with $A$-valued pre-inner product $\lprod A {\cdot} {\cdot} : M \times M
\mapsto A$, linear in first argument and conjugate linear in second, such that
$\lprod A {a x} y = a \lprod A x y $ for $x,y \in M, a \in A$ and the
two-sided ideal $\Img \lprod A {\cdot} {\cdot}$ is dense in $A$.
\end{definition}

Note that $\lprod A {a x} y = a \lprod A x y $ for an $A$-valued inner product
imply also $\lprod A x {a y} = \lprod A x y a^*$. That is why $\Img \lprod A
{\cdot} {\cdot}$ is two-sided ideal as mentioned. Definition of a right rigged
space we obtain by simple reflection from the left to the right:

\begin{definition}
A right $A$-module $M$ is called a right $A$-rigged space if it is
endowed with pre-inner product $\rprod A {\cdot} {\cdot} : M \times M \mapsto
A$,
conjugate linear in first argument and linear in second, such that $\rprod A x
{ya} = \rprod A x y a$ for $x,y \in M, a \in A$ and the ideal $\Img \rprod A
{\cdot} {\cdot}$ is dense in $A$.
\end{definition}

Let now $A,B$ be pre-$\Cst$-algebras.
\begin{definition}\label{imprimitivity}
An $A-B$ bimodule $M$ is called an imprimitivity bimodule if

(1) $M$ is left-$A$-right-$B$-rigged space;

(2) $\lprod A x y z = x \rprod B y z$;

(3) $\rprod B {ax}{ax} \le \norm{a}^2_A \; \rprod B x x$ and $\lprod A {xb}{xb}
\le \norm{b}^2_B \; \lprod A x x$

for $x,y,z \in M$ and $a \in A$, $b \in B$.
\end{definition}

Note that in imprimitivity bimodule we also have relation $\lprod A x {yb} =
\lprod A {xb^*} y$ for any $x,y \in M$ and $b \in B$.
Indeed, for $b \in \Img \rprod B {\cdot}{\cdot}$ it is a consequence of 
relation~(2) in definition above. Let us check necessary
continuity. Suppose $\norm {b_n} \rightarrow 0$. Then $\norm {\lprod A {y b_n}{y
b_n}} \le \norm{b_n}^2 \norm{\lprod A y y}
\rightarrow 0$ by (3). Now by Proposition~2.9 in \cite{Riff1974} we have
$\norm{\lprod A x {yb_n}} \le \norm{\lprod A x x}^{\frac 12} \norm{\lprod A
{yb_n} {yb_n}}^{\frac 12}
\rightarrow 0$, and analogously $\norm{\lprod A {x b_n^*} y} \rightarrow 0$. 
Evidently, symmetrical relation
$\rprod B {ax} y = \rprod B x {a^*y}$ for any $x,y \in M$ and $a \in A$ also
holds.

\begin{definition}{ \emph (\cite{Riff1976}) }
Two pre-$\Cst$-algebras $A$, $B$ are said to be strongly Morita equivalent if
there exist an $A-B$-imprimitivity bimodule.  
\end{definition}

\begin{example}\label{free_mod}
Let $A$ be a pre-$\Cst$-algebra with $1$. Consider $M = A^n$~--- free right $A$-
module of rank $n$. Then $End_A M = \Mat n A$~--- the ring of $n \times n$
matrices with entries in $A$, which is a unital pre-$\Cst$-algebra again. 
Then $M$ is $\Mat n A - A$-imprimitivity bimodule with inner
products
$$
\rprod A x y = x^* y = \sum_{i} {x_i^*}y_i
$$
$$
\lprod {\Mat n A} x y = x y^* = (x_i y_j^*)_{i,j=1}^n
$$
\end{example}

\begin{example}\label{subgroups}
Let $G$ be a locally compact group, and let $H$ and $K$ be closed subgroups 
of $G$. Let $A = \Cst(K, G \slash H)$, $B = \Cst(H, K \backslash G)$ be
transformation group $\Cst$-algebras for left action of $K$ on $G\slash H$
and right action of $H$ on $K \backslash G$ correspondingly. It is shown in
\cite{Riff1976} that there is a natural $A-B$-imprimitivity bimodule $M$~--- a
completion of the space $C_c(G)$ of $\Co$-valued continuous functions with
compact support on $G$ with respect to an appropriate norm, with inner
products given on $f,g \in C_c(G)$ by:
$$
\lprod A f g (k,x) = \beta(k)\int_H f(\tilde{x}h)g^*(h^{-1}\tilde{x}^{-1}k)dh
$$
where $\tilde{x} \in G$ is any representative of class $x$, i.e. $x=\tilde{x}H$,
$$
\rprod B f g (h,y)= \gamma(h)\int_K f^*(\tilde{y}^-1 k)g(k^{-1}\tilde{y}h)dk 
$$ 
where $y = K \tilde{y}$. Here $\beta(\cdot)= \left( \frac{\delta_G(\cdot)}
{\delta_K(\cdot)} \right)^{\frac 12}$, $\gamma(\cdot) =
\left(\frac{\delta_G(\cdot)}{\delta_H(\cdot)} \right)^{\frac 12}$, $\delta_{G}$,
$\delta_{H}$, $\delta_{K}$ are the modular functions of locally compact groups
$G,H,K$ correspondingly, an involution is defined on $C_c(G)$ by $g \mapsto
g^*(z) = \delta_G(z^{-1})\bar{g}(z^{-1})$, and all integrals above are taken
w.r.t. left Haar measures. 
\end{example}

We will see in Sections 2,3 below that this two examples are quite similar.

Strong Morita equivalence implies Morita equivalence, i.e. equivalence of
categories of Hermitian representations(\cite{Riff1974}). It is not obvious
from definitions that strong Morita equivalence is indeed an equivalence
relation. In \cite{Riff1974} an inverse imprimitivity bimodule is constructed,
showing that this relation is symmetric. In Section~4 we define a natural
structure of imprimitivity bimodule on tensor product of imprimitivity
bimodules for unital $\Cst$-algebras. In particular it makes evident
transitivity of strong Morita equivalence for unital $\Cst$-algebras.

\section{Inner products for projective module}

We generalize Example~\ref{free_mod} in current section. On a projective
module pre-inner products which satisfy all algebraic
relations from Definition~\ref{imprimitivity} were introduced
in~\cite{AG0202109}. We are going to check the condition of density of images
for these inner products now. 

Let $A$ be $\Cst$-algebra with $1$, let $p \in \Mat n A$ be projection, i.e.
$p=p^*=p^2$. Consider a submodule $M = p A^{n}$ of right $A$-module $A^n$
consisting of such columns which are invariant under left multiplication 
by $p$. Then $End_A M = p \Mat n A p$, where matrices act by multiplication
from the left. $p \Mat n A p$ is a $\Cst$-algebra with norm restricted from
$\Mat n A$, and since $\norm{p}=1$ this is a unital $\Cst$-algebra. We consider
two inner products on $M$ which are restrictions of inner products from
Example~\ref{free_mod}:
$$
\rprod A x y = x^* y = \sum_{i} {x_i^*}y_i 
$$
$$
\lprod {p \Mat n A  p} x y = x y^* = (x_i y_j^*)_{i,j=1}^n 
$$
Then $\Img \lprod {p \mat A p} {\cdot}{\cdot} = p \Mat n A p$, and
$\Img \rprod A {\cdot}{\cdot} = \sum_{i,j} A p_{i,j}A$~--- the ideal in $A$
generated by matrix entries of $p$. 

\begin{proposition}
$p A^n$ with inner products defined above is $p \Mat n A p -
A$-imprimitivity bimodule if and only if $\Img \rprod A {\cdot}{\cdot} = A$.
\end{proposition}
\begin{proof}
In unital $\Cst$-algebra $A$ there is no dense ideal except $A$. So, condition 
$\Img \rprod A {\cdot}{\cdot} = A$ is necessary for $p A^n$ to be
right $A$-rigged space. We show it is sufficient. Indeed, all necessary
identities
for $\rprod A {\cdot}{\cdot}$ and $\lprod {p \mat n A p} {\cdot}{\cdot}$ are
satisfied as they are satisfied in Example~\ref{free_mod}, and we already
mentioned that $\Img \lprod {p \mat n A p} {\cdot}{\cdot} = p \Mat n A p$.
\end{proof}

In the next section we show that any imprimitivity bimodule between unital
$\Cst$-algebras is of this form. This idea also comes from works of M.~
Rieffel~--- one may compare Theorem~\ref{general_bimodule} below to
Proposition~2.1 in \cite{Riff1981}. 

\section{Imprimitivity bimodule for $\Cst$-algebras with 1}

\begin{theorem}\label{general_bimodule}
Let $A, B$ be two strongly Morita equivalent $\Cst$-algebras with 1, and $M$ be
a $B-A$-imprimitivity bimodule. Then 

(1) $B = End_A M$; 

(2) there exist $n \in \Z$, projection $p \in \Mat n A$ and isomorphism of
right $A$-modules $\Psi: M \rightarrow p A^{\infty}$ such that for $u,v \in
M$: 
$$
\rprod A u v = \Psi(u)^*\Psi(v)
$$ 
$$
\lprod B u v = \Psi^{-1} \circ \Psi(u)\Psi(v)^* \circ \Psi
$$
\end{theorem}
\begin{proof}
As $B$ is unital $\Cst$-algebra, any dense ideal in it is $B$. Then there exist
integer $n$ and $x_1,\dots,x_n,y_1,\dots,y_n \in M$ such that
$$
1_B = \sum_i \;\lprod B {x_i}{y_i}.
$$ 

Consider unital $\Cst$-algebra $C = \Mat n A$ and $B-C$-bimodule $N = M^n$
consisting of columns of elements of $M$. We define inner products on $N$ by
$$
\lprod B m n = \sum_i \; \lprod B {m_i} {n_i} 
$$
$$
\rprod C m n = (\rprod A {m_i}{n_j})_{i,j=1}^n 
$$
One can check $N$ is $B-C$-imprimitivity bimodule. Now $\lprod B x y = 1_B$.
Consider $z = x \rprod C y y ^{1/2}$. Then 
$$
\lprod B z z = \lprod B x {x \rprod C y y} = \lprod B x {\lprod B x y y}
= \lprod B y x \lprod B x y = 1_B,
$$
and $p = \rprod C z z$ is a projection. Indeed, obviously $p^*=p$ and 
$$
\rprod C z z \rprod C z z = \rprod C z {z \rprod C z z} = \rprod C z {\lprod B
z z z} = \rprod C z z.
$$

Consider homomorphism of right $A$-modules $\Psi: M \rightarrow p A^n$,
$\Psi(m) = (\rprod A {z_i} m)$, and unital homomorphism of rings 
$\Phi: B \rightarrow p \Mat n A p$, $\Phi(b) = \rprod C z {bz}$. We now prove
they are both correctly defined and are in fact isomorphisms.

For $\Psi$ consider $j: M \rightarrow N$ given by $j(m) = (m
\delta_{1i})_{i=1}^n$. Then columns of $\rprod C z {j(m)}$ are
$\Psi(m),0,\dots,0$. Since $P \rprod C z {j(m)}$ 
$= \rprod C {z \rprod C z z} {j(m)}$\\$=\rprod C {\rprod C z z z} {j(m)}
=\rprod C z {j(m)}$, $\Psi(m) \in p A^n$. Injectivity of $\Psi$ follows 
from $z \rprod C z n = \lprod B z z n = n \ne 0$ for nonzero $n \in N$.
Surjectivity follows from fact that $\Img \rprod C z {\cdot} = p C$.

$\Phi(b)$ is obviously invariant under left ant right multiplication by $p$.
To prove injectivity we note that $b = \lprod B y x b \lprod B x y = \lprod B
{bz} z$, so $b z \ne 0$ if $b \ne 0$. Thus also $\rprod C z {bz} \ne 0$.
Surjectivity follows from fact that $p C p$ is spanned by $p \rprod C m n p$ and
equality $p \rprod C m n p =\Phi(\lprod B m z \lprod B n z)$. Also 
$\Phi(b_1b_2) = \Phi(b_1) \Phi(b_2)$.

To prove statement it remains to check that $\rprod A  u v = \Psi(u)^*\Psi(v)$,
$\Phi(\lprod B u v) = \Psi(u)\Psi(v)^*$ and $\Phi(\lprod B u v) \Psi(t) =
\Psi(\lprod B u v t)$. Indeed, 
$$\Psi(u)^*\Psi(v) = \sum_{i} \; \rprod A u {z_i}
\rprod A {z_i} v = \rprod A u {\lprod B z z v} = \rprod A u v .$$
Next, we compare $(i,j)$'th matrix entry for $\Phi(\lprod B u v) = \rprod C z
{\lprod B u v z}$ and $\Psi(u)\Psi(v)^*$:
$$
\rprod A {z_i} {\lprod B u v z_j} = \rprod A {z_i} u \rprod A v {z_j}.
$$
Now we compare $i$'th coordinate in $\Phi(\lprod B u v) \Psi(t)$ and
$\Psi(\lprod B u v t)$:
$$
(\Psi(u)\Psi(v)^*\Psi(t))_i = \Psi(u)_i \rprod A v t
= \rprod A {z_i} u \rprod A v t = \rprod A {z_i} {\lprod B u v t}.
$$
\end{proof}

\begin{corollary}\label{one_defines_another}
Suppose there are two structures of a $B-A$-imprimitivity bimodule on bimodule
$M$: $\lprod B{\cdot}{\cdot} ^i$ and $\rprod A{\cdot}{\cdot} ^i$ for $i = 1,2$.
If $\rprod A{\cdot}{\cdot} ^1 = \rprod A{\cdot}{\cdot} ^2$ then also\\
$\lprod B{\cdot}{\cdot} ^1=\lprod B{\cdot}{\cdot} ^2$, and vice
versa.
\end{corollary}
\begin{proof}
Due to theorem above it is sufficient to check the statement in case $M = p
A^{\infty}$ and $\rprod A x y ^1 = \rprod A x y ^2 = x^*y$. Then for any $z \in
p
A^{\infty}$ we have $\lprod B x y z = x \rprod A y z = x y^* z$. Taking $z =
p_k$ for all columns of $p = (p_k)$ we get $\lprod B x y = \lprod B x y p =
x y^* p = x y^*$, so the second inner product is defined by the first one. 
\end{proof}

\section{Composition of strong Morita morphisms}

Evidently choice of inner products for an $A-B$-imprimitivity bimodule $M$ is
non-unique. For example, we can multiply them both by positive number and with 
such new inner products $M$ will be again an $A-B$-imprimitivity bimodule.
Anyway, following theorem gives one natural choice of imprimitivity bimodule
structure on tensor product of two imprimitivity bimodules. 

\begin{theorem}\label{tensor}
Let $A,B,C$ be unital $\Cst$-algebras, $M,N$ be $A-B$ and $B-C$-imprimitivity
bimodules correspondingly. Then $M \underset{B} \otimes N$ with inner
products defined by
$$
\rprod C {x\otimes z} {y\otimes t} = \rprod C z {\rprod B x y t} 
$$
$$
\lprod A {x\otimes z} {y\otimes t} = \lprod A {x \lprod B z t} y 
$$
is an $A-C$-imprimitivity bimodule.
\end{theorem}
\begin{proof}
Let $K = M \underset{B} \otimes N$. We check that $K$ is a right $C$-rigged
space. First, let us see that $\rprod C {\cdot} {\cdot}$
on $K$ is well-defined  $C$-valued inner product antilinear in first variable.
Indeed, for $b \in B$
$$
\rprod C {xb \otimes z}{y \otimes t} = \rprod C z {\rprod B {xb}y t} = \rprod C
z {b^* \rprod B x y t} 
$$
$$
= \rprod C {bz} {\rprod B x y t} = \rprod C {x \otimes bz}{y \otimes t}.
$$
Taking $b \in \Co 1$ we see that $\rprod C {\cdot}{\cdot}$ is antilinear in
first variable. Analogously 
$\rprod C {x \otimes z}{yb \otimes t} = \rprod C {x \otimes z}{y \otimes b t}$
and $\rprod C {\cdot}{\cdot}$ is linear 
in second variable. To see positivity of $\rprod C {\sum_{i=1}^n x_i \otimes
z_i}{\sum_i x_i \otimes z_i} = 
\sum_{i,j} \rprod C{z_i}{\rprod B {x_i}{x_j}{z_j}}$, we recall that $M^n$ is an
$A - \Mat n B$-imprimitivity 
bimodule. So matrix $H = \rprod B {x_i}{x_j}$ is positive element of $\Mat n B$.
 Thus $\rprod C {y} {H y} 
= \sum_{i ,j}\rprod C {y_i}{h_{i,j}y_j} \ge 0$ as $N^n$ is an $\Mat n B -
C$-imprimitivity bimodule.  

Consider elements $x_i, y_i$ in $M$ such that $\sum_i \rprod B {x_i} {y_i} = 1$.
Then for any $z,t \in N$
$$
\sum_i \rprod C {x_i \otimes z}{y_i \otimes t} = \rprod C z t, 
$$
so $\Img \rprod C {\cdot}{\cdot}$ on $N$ is a subset of $\Img \rprod C
{\cdot}{\cdot}$ on $K$. Thus 
$\Img \rprod C {\cdot}{\cdot}$ on $K$ is dense in $C$.

For $c \in C$ we obviously have relation 
$$
\rprod C {x \otimes z}{y \otimes t c} = \rprod C {z}{\rprod B x y  t c} = \rprod
C {x \otimes z}{y \otimes t} c,
$$
so we proved $K$ is a right $C$-rigged space. Analogously $K$ is a left
$A$-rigged space.

Now we check condition (2) in Definition~\ref{imprimitivity} :
$$
v \otimes w \rprod C {x \otimes z}{y \otimes t} = v \otimes w \rprod C {\rprod B
y x z}{t}=
v \otimes \lprod B w {\rprod B y x z} t
$$
$$
= v \lprod B w {\rprod B y x z} \otimes t = v \lprod B w z \rprod B x y \otimes
t = 
\lprod A {v \lprod B w z} x y \otimes t
$$ 
$$
= \lprod A {v \otimes w}{x \otimes z} y \otimes t.
$$

For condition (3) consider $a \in A$ and
$$
\rprod C {a \sum_{i=1}^n {x_i} \otimes {z_i}}{a \sum_i {x_i} \otimes {z_i}} =
\sum_{i,j} \rprod C {z_i} {\rprod B {ax_i}{ax_j} z_j}
$$
$$
\le \sum_{i,j} \rprod C {z_i} {\norm{a}^2 \rprod B {x_i}{x_j} z_j} = \norm{a}^2
\rprod C {\sum_i {x_i} \otimes {z_i}}{\sum_i {x_i} \otimes {z_i}} 
$$
as $M^n$ is an $A - \Mat n B$-imprimitivity bimodule. Analogously we can check
condition (3) for $\lprod A {\cdot}{\cdot}$.
\end{proof}

We remark that this statement is also true in case of unital
pre-$\Cst$-algebras. The proof is the same just we need additional 
continuity arguments to prove density of images of  pre-inner products.

\section{Morita bimodules over quantum tori}

Recall that quantum torus $A_{\theta}$ for $\theta \in \R \backslash \Q$ is 
a transformation group $\Cst$-algebra $\Cst(\theta \Z,\R \slash \Z)$. It is
known that $A_{\theta}$ is universal $\Cst$-algebra generated by two unitaries
$U,V \in
A_{\theta}$ satisfying relation $UV = e(\theta)VU$. The choice of such
unitaries is not unique. If $U,V \in A_{\theta}$ are chosen we call them a
frame.

From Example~\ref{subgroups} we see that $A_{\theta} = \Cst(\theta \Z,\R \slash
\Z)$ is strongly Morita equivalent to $\Cst(\Z,\R \slash \theta \Z) \cong
\Cst(\frac 1{\theta}\Z,\R \slash \Z) = A_{\frac 1 {\theta}}$. Obviously
$A_{\theta + 1} = A_{\theta}$, as relation $UV = e(\theta)VU$ is invariant
under transformation $\theta \mapsto \theta + 1$. Also $A_{\theta} \cong
A_{-\theta}$ as we can map $U$ to $V^{'}$ and $V$ to $U'$ for any frames
$U,V \in A_{\theta}$, $U^{'},V^{'} \in A_{\theta}$. Recall that $GL_2(\Z)$ acts
on complex numbers by 
$$
\begin{pmatrix}
a&b\\c&d\\
\end{pmatrix}
\theta = \frac{a \theta + b}{c \theta + d}.
$$
So we see that $A_{\theta}$ is strongly Morita equivalent to $A_{g \theta}$ for
any $g \in GL_2(\Z)$.
Indeed, as $GL_2(\Z)$ is generated by $\begin{pmatrix}0&1\\1&0\\ \end{pmatrix}$
and $\begin{pmatrix}1&1\\0&1\\ \end{pmatrix}$, its orbit is generated by
transformations $\theta \mapsto \theta+1$ and $\theta \mapsto \frac 1{\theta}$.
Conversely, it is shown in \cite{Riff1981} that $A_{\theta}$ and
$A_{\theta^{'}}$ are not strongly Morita equivalent if $\theta$ and $\theta^{'}$
don't lie in the same orbit of $GL_2(\Z)$. 

Below we recall an explicit construction of
$A_{g\theta}-A_{\theta}$-imprimitivity
bimodule $E(g,{\theta})$ for $g \in SL_2(\Z)$, $\theta \in \R \backslash \Q$
(\cite{QA0203160}, \cite{QA0211262},\cite{AG0202109}). (Bimodules for $g \in
GL_2(\Z)$ can be easily obtained from those by composition with homomorphism $U
\mapsto V^{'}$,$V\mapsto U^{'}$ of quantum tori on the left.) It is proven in
\cite{QA0203160} that $E(hg,\theta) \cong E(h, g \theta) \underset{A_{g
\theta}}\otimes E(g, \theta)$ as bimodule, and we claim in theorem below that
inner products satisfy relations of Theorem~\ref{tensor}. 

To construct our bimodules we need to fix a frame in $A_{\theta}$ for each
$\theta \in \R \backslash \Q$. (If $\theta = \theta^{'}$ modulo $\Z$ then
the frames should coincide.) Let $g = \begin{pmatrix}a&b\\c&d\\\end{pmatrix} \in
SL_2(\Z)$. If $c = 0$ we put $E(g,\theta)=A_{\theta}$ with action of $A_{g
\theta} = A_{\theta}$ via multiplication from the left and action of
$A_{\theta}$ by right multiplication, and define inner products by $\lprod {A_{g
\theta}} a b = a b^*$ and $\rprod {A_{\theta}} a b = a^* b$ as in
Example~\ref{free_mod}. If $c \ne 0$, we consider the space $E^0(g,\theta)=S(\R
\times \Z\slash c\Z)$ with following actions of generators $U,V$ of $A_{\theta}$
and $U^{'}, V^{'}$ of $A_{g \theta}$ on $f \in E^0(g,\theta)$:
$$
(f U)(x,\alpha) = f(x-\frac{c \theta + d}c,\alpha - 1)
$$
$$
(f V)(x, \alpha) = e(x - \alpha \frac d c )f(x,\alpha) 
$$
$$
(U^{'} f)(x,\alpha) = f(x-\frac 1 c,\alpha - a)
$$
$$
(V^{'} f)(x, \alpha) = e(\frac x {c \theta + d} - \frac {\alpha} c )f(x,\alpha) 
$$
We define for $f,s \in E^0(g,\theta)$ inner products:
$$
\lprod {A_{g \theta}} f s = \sum_{n \in \Z^2} \iprod f
{U^{'n_1}V^{'n_2}s}_{L_2} U^{'n_1}V^{'n_2} 
$$  
$$
\rprod {A_{\theta}} f s = \frac 1 {c \theta + d} \sum_{n \in \Z^2} \iprod s
{f U^{n_1}V^{n_2}}_{L_2} U^{n_1}V^{n_2} 
$$  
Let $E(g,\theta)$ be the completion of $E^{0}(g,\theta)$ with the norm $\norm f
= \norm{\lprod {A_{g\theta}} f f}^{\frac 1 2}$. Then $E(g, \theta)$ is an $A_{g
\theta}-A_{\theta}$-imprimitivity bimodule (Theorem 3.2 in \cite{AG0202109}).

In \cite{QA0203160},\cite{QA0211262} there are constructed bimodule isomorphisms
$t_{h,g}: E(h,g\theta) \underset{A_{g \theta}}\otimes E(g, \theta) \rightarrow
E(hg,\theta)$ for $h,g \in SL_2(\Z)$. 

\begin{theorem}\label{tensor_for_tori}
For $h,g \in SL_2(\Z)$, $f_1, s_1 \in E(h,g \theta)$ and 
$f_2, s_2 \in E(g, \theta)$
$$
\lprod {A_{hg \theta}} {t_{h,g}(f_1 \otimes f_2)} {t_{h,g}(s_1 \otimes s_2)} =
\lprod {A_{hg \theta}} {f_1 \; \lprod {A_{g \theta}} {f_2} {s_2}} {s_1}
$$
$$
\rprod {A_{\theta}} {t_{h,g}(f_1 \otimes f_2)} {t_{h,g}(s_1 \otimes s_2)} =
\rprod {A_{\theta}} {f_2} { \rprod {A_{g \theta}} {f_1} {s_1} {s_2}}
$$
\end{theorem}
\begin{proof}
First, due to Theorem~\ref{tensor} and Corollary~\ref{one_defines_another} it is
enough to check only one of two statements of the theorem. We prefer the 
second one. 

As maps $t_{h,g}$ are associative (Proposition 1.2 in \cite{QA0211262}) 
it is enough to check the statement for generators of $SL_2(\Z)$ at
place of $h$ only.  Indeed, suppose the statement is true for
$E(h_1,g\theta)\otimes E(g,\theta)$, $E(h_2,h_1g\theta)\otimes E(h_1g,\theta)$
and $E(h_2,h_1g\theta)\otimes E(h_1,g\theta)$. Then it is true for
$E(h_2h_1,g\theta)\otimes E(g,\theta)$ due to associativity relation
$$
t_{h_2h_1,g} \circ (t_{h_2,h_1} \otimes id) =
t_{h_2,h_1g} \circ (id \otimes t_{h_1,g}). 
$$

Take $h = \begin{pmatrix}1&1\\0&1\\ \end{pmatrix}$. Then $f_1, s_1 \in A_{g
\theta}$, $\rprod {A_{g \theta}} {f_1} {s_1} = f_1^* s_1$, $t_{h,g}(f_1, f_2) =
f_1 f_2$ (in sense of left action) and similar $t_{h,g}(s_1,s_2) = s_1 s_2$.
As $h \begin{pmatrix}a&b\\c&d\\\end{pmatrix}= \begin{pmatrix}
a+c&b+d\\c&d\\\end{pmatrix}$ we have no changes in formulas for
action of quantum tori, so $E(hg,\theta) = E(g,\theta)$ and
$$
\rprod {A_{\theta}} {f_1 f_2}{s_1 s_2} = \rprod {A_{\theta}} {f_2}{f_1^*s_1
s_2} 
$$
as $E(g,\theta)$ is an $A_{g \theta}-A_{\theta}$-imprimitivity bimodule. Indeed,
for an $A-B$-imprimitivity bimodule $M$ we have $\rprod B {ax} y =
\rprod B x {a^* y}$ for $a \in A$, $x,y \in M$.

Now take $h = \begin{pmatrix}0&-1\\1&0\\ \end{pmatrix}$. Then $hg =
\begin{pmatrix}-c&-d\\a&b\\ \end{pmatrix}$. Let us consider the case $g \ne
h$, $c \ne 0$. Cases $g = h$ and $c = 0$ can be done analogously. Obviously we
can restrict to dense set of Schwartz functions $f_1, s_1 \in
E^0(h,g \theta)$, $f_2, s_2 \in E^0(g, \theta)$. $E^0(h,g \theta) = S(\R)$ with
$$
\rprod {A_{g \theta}} {f_1} {s_1} = \frac 1{g \theta} \sum_{n \in \Z^2} \int
s_1(y)e(-y n_2)\bar{f_1}(y - n_1 \theta) dy \; U^{n_1}V^{n_2}
$$
where $U,V \in A_{g \theta}$. Let $U^{'}$,$V^{'}$ be generators of $A_{\theta}$.
Comparing coefficients near $U^{'m_1}V^{'m_2}$ in identity, which we need to
prove, we see that it is equivalent to 
$$
\frac 1 {a \theta + b} \iprod {t_{h,g}(s_1 \otimes s_2)}{t_{h,g}(f_1 \otimes f_2
U^{'m_1}V^{'m_2})}_{L_2} = \frac 1 {c \theta + d} \iprod { \rprod {A_{g \theta}}
{f_1}{s_1} s_2}{f_2 U^{'m_1}V^{'m_2}}_{L_2} 
$$
Substituting $f_2$ instead of $f_2 U^{'m_1}V^{'m_2}$, we need to prove for
arbitrary $f_1, s_1 \in S(\R)$, $f_2, s_2 \in S(\R \times \Z \slash c \Z)$
$$
\iprod {t_{h,g}(s_1 \otimes s_2)}{t_{h,g}(f_1 \otimes f_2)}_{L_2} = 
\sum_{n \in \Z^2}\int s_1(y)e(-y n_2)\bar{f_1}(y - n_1 \theta) dy \;
\iprod {U^{n_1}V^{n_2} s_2}{f_2}_{L_2}
$$ 
This is a routine computation using Poisson summation formula. We use
abbreviations LHS (RHS) for left-(right)-hand side of this identity
correspondingly. By explicit formula for $t_{h,g}$ (Proposition 1.2 in
\cite{QA0211262})
$$
t_{h,g}(s_1 \otimes s_2)(x,\alpha) = \sum_{n \in \Z} s_1\left( \frac x{c \theta
+ d} + g\theta \left( \frac{cb}{a} \alpha - n\right) \right) s_2\left( x 
- \frac b a \alpha + \frac n c, a n\right),
$$
and analogously for $t_{h,g}(f_1 \otimes f_2)$. Now
$$
LHS = \sum_{n,m \in \Z} \sum_{\alpha \in \Z \slash a \Z} \int s_1(z)s_2(y -
\frac{m-n}c, an) \bar{f_1}(z - g\theta(m-n))\bar{f_2}(y,a m) dy
$$
where $z = \frac x{c \theta + d} + g\theta \left( \frac{cb}{a} \alpha -
n\right)$ and $y = x - \frac b a \alpha + \frac m c$. Let us represent $m
= d m_1 + c m_2$ with $m_1 \in \Z \slash c \Z$ and $m_2 \in \Z$. Then $am = m_1$
and $an = m_1 - a(m-n)$ modulo $c$. Introducing new variable $n_1 = m-n$ we
proceed:
$$
= \sum_{m_1 \in \Z \slash c \Z} \int \sum_{n_1 \in \Z} \sum_{m_2 \in \Z, \alpha
\in \Z \slash a \Z} s_1(z)\bar{f_1}(z - g\theta n_1)
(U^{n_1}s_2)(y,m_1)\bar{f_2}(y,m_1)dy
$$
Let us express $z$ via $y$ and summing variables: 
$$
z = \frac 1{c \theta+d}\left( y + \frac ba \alpha - \frac m c \right) + g
\theta(\frac{cb}a \alpha - n)
$$
$$
= \frac 1{c \theta+d}\left( y + \frac ba \alpha - m_2 - \frac d c m_1 \right) +
\frac{a \theta + b}{c \theta + d}(\frac{cb}a \alpha - c m_2 - d m_1 + n_1)
$$
$$
= (b \alpha - m_2 a)-\frac{ad}c m_1 + \frac 1{c\theta+d}(y + (a \theta + b)n_1)
$$
Denote $n_2 = b \alpha - m_2 a$, and $z_0 = z - n_2$. Then by Poisson
summation formula
$$
\sum_{n_2 \in \Z} s_1(n_2 + z_0)\bar{f_1}(n_2+z_0- g \theta n_1) = 
\sum_{n_2  \in \Z} e(z_0)^{n_2}\int e(-t n_2)s_1(t)\bar{f_1}(t - g \theta n_1)
dt.
$$
We put this into LHS, and note that $e(z_0)^{n_2}(U^{n_1}s_2)(y,m_1) =
(U^{n_1}V^{n_2}s_2)(y,m_1)$. So LHS =
$$
\sum_{m_1 \in \Z \slash c\Z} \int \sum_{n_1,n_2 \in \Z} \int
e(-t n_2)s_1(t)\bar{f_1}(t - g \theta n_1)dt
(U^{n_1}V^{n_2}s_2)(y,m_1)\bar{f_2}(y,m_1)dy
$$
$$
= \sum_{n_1,n_2} \int e(-t n_2)s_1(t)\bar{f_1}(t - g \theta n_1)dt \iprod
{U^{n_1}V^{n_2}s_2}{f_2}_{L_2} = RHS
$$
\end{proof}

\section{Real multiplication}

Irrational number $\theta \in \R \backslash \Q$ is a root of quadratic equation
if and only if there exist matrix $g \in SL_2(\Z)$, $g \ne \pm 1$ such that $g
\theta = \theta$. Let us fix such $g$ and $\theta$. It follows from Section~5
that there are nontrivial $A_{\theta}-A_{\theta}$-imprimitivity bimodules
exactly in this case. Now we are going to construct a graded ring
$R=R(g,\theta)=\underset{n \ge 1}\oplus R_n$ using tensor products and inner
products in these imprimitivity bimodules. We start with construction of another
graded ring due to Polishchuk \cite{AG0212306}, which uses only tensor products.
 
We consider the set of bimodules $E(g^n,\theta)$, $n \ge 1$ defined in
previous section, and have the family of isomorphisms 
$$
t_{g^m,g^n} : E(g^n,\theta) \otimes E(g^m,\theta) \widetilde{\rightarrow}
E(g^{n+m},\theta).
$$ 
Let $\Siegel k = \{ M \in \Mat k {\Co} | M=M^t \;\text{and}\; \Im(M)>0 \}$ be
so-called Siegel upper half-plane. So, $\Siegel 1$ is just an upper half of
complex
plane $\Co$, and we fix $\tau \in \Siegel 1$. Denote matrix entries of $g^n$ by
$\begin{pmatrix}a_n&b_n\\c_n&d_n\\ \end{pmatrix}$. Denote $\mu_n = \tau
\frac{c_n}{c_n \theta + d_n}$.
Note that $c_n \theta + d_n$ is an eigenvalue of $g^n$, so it is nonzero. Also
$c_n \ne 0$ as $g^n$ is a nontrivial matrix 
stabilizing $\theta$. Thus $\mu_n \ne 0$. Denote 
$$
E_n = \begin{cases}
\left\{ \phi_f(x,\alpha) = e(\mu_n \frac{x^2}2)f(\alpha) \Big| f :
\Z \slash c_n \Z \rightarrow \Co \right\}, &\frac{c_n}{c_n \theta + d_n} > 0\\
\left\{ 0 \right\}, &\frac{c_n}{c_n \theta + d_n} < 0\\
\end{cases}
\subset E(g^n,\theta).
$$
$E_n$ is either $0$ or a $|c_n|$-dimensional vector space. In fact we have
either $E_n = \{ 0 \}$ for all $n$
or $E_n \ne \{0 \}$ for all $n$. Indeed, we see that definition of $E_n$ is the
same for $E(g,\theta)$ and $E(-g,\theta)$.
Thus taking either $g$ or $-g$ instead of $g$ we can suppose that $c_1 \theta +
d_1 > 0$. $c_1 \theta + d_1$ is an 
eigenvalue of $g$, so $g$ has positive eigenvalues. Now it follows from
$\sum_{n=1}^{\infty}c_n t^n = \frac{c t}{t^2 - tr(g) t + 1}$
that all $c_n$ have the same sign, as all coefficients of power series for 
$\frac 1{t^2 - tr(g) t + 1}$ are positive. All $c_n \theta + d_n$
are eigenvalues of $g^n$, so they are also positive.

Consider the set
$$
S_{\theta} := \left\{ g=\begin{pmatrix}a&b\\c&d\\ \end{pmatrix} \in  SL_2(\Z)
\Big | g \ne \pm 1, g \theta = \theta,
tr(g)>0 \; \& \; c>0  \right \}
$$
It is always nonempty: we already showed how to satisfy first three
conditions, then if the forth is not satisfied we can take $g^{-1}$ instead
of $g$. 

Further we suppose $g \in S_{\theta}$. Then all $E_n$ are nonzero vector spaces.
It was noticed already in \cite{QA0203160}
that vector spaces $E_n$ are preserved under tensor products of bimodules.
Following can be checked by direct computation:
\begin{proposition} For $f : \Z \slash c_n \Z \rightarrow \Co$, $g:\Z \slash
c_m \Z \rightarrow \Co$ we have $t_{g^n,g^m}(\phi_f \otimes \phi_g) = \phi_{f
\underset{n,m}\star g}$ where
$$
f \underset{n,m}\star g (\alpha) = \sum_{q \in \Z}e\left(\frac{\tau}2
\frac{c_{n+m}}{c_n c_m}\left( q - \frac{c_m d_{n+m}}{c_{n+m}} \alpha
\right)^2\right)f(a_nd_{n+m}\alpha -q)g(a_m q) 
$$
is a function on $\Z \slash c_{n+m} \Z$.
\end{proposition}

Now we consider the graded ring $E = \oplus_{n \ge 1} E_n$ with multiplication
given by $\phi_f * \phi_g := \phi_{t \underset{n,m}\star g} \in E_{n+m}$ for
$\phi_f \in E_n$, $\phi_g \in E_m$. Associativity of this multiplication follows
from identity
$$
t_{g^{n+m},g^k} \circ (t_{g^n,g^m} \otimes id) =
t_{g^{n},g^{m+k}} \circ (id \otimes t_{g^m,g^k}) : E_n
\otimes E_m \otimes E_k \rightarrow E_{n+m+k}
$$
stated in Proposition 1.2 in \cite{QA0211262}. Note that if we choose for basis
in $E_n$ functions of the form $\phi_f$ with characters $f \in (\Z \slash c_n
\Z)^*$, we would get multiplication table consisting of values at rational
points of various theta functions with rational characters
$\theta\left[\begin{aligned}
\alpha\\\beta\\\end{aligned}\right](\gamma,\delta \tau)$ where
$\alpha,\beta,\gamma,\delta \in \Q$ (see, e.g. \cite{Mumf1983}). For example,
$$
1 \underset{n,m}\star 1 (\alpha) = \theta\left[
\begin{aligned} \frac{c_m d_{n+m}}{c_{n+m}}\alpha \\0\\ \end{aligned}
\right]\left(0,\frac{\tau}2 \frac{c_{n+m}}{c_n c_m}\right).
$$

In \cite{AG0212306} (Theorem 2.4) there are established criterions whether the
ring $E$ is generated over $\Co$ by $E_1$, is quadratic and is Koszul. Using
them we state the criterion
whether there exist $g \in S_{\theta}$ such that $E$ have these good properties:

\begin{theorem}\label{good_E}
Let $\theta \in \R \backslash \Q$ be a quadratic irrationality, and $\theta^{'}$
be its Galois conjugate. Then the following conditions are equivalent:\\
(1) $|\theta - \theta^{'}|<1$;\\
(2) there exist $g \in S_{\theta}$ such that the ring $E$ is generated by $E_1$
over $\Co$;\\
(3) there exist $g \in S_{\theta}$ such that the ring $E$ is quadratic;\\
(4) there exist $g \in S_{\theta}$ such that the ring $E$ is Koszul.\\
\end{theorem}
\begin{proof}
First we show (2),(3) and (4) imply (1). Let $g =
\begin{pmatrix}a&b\\c&d\\ \end{pmatrix}$ with given properties exist. As $g \in
S_{\theta}$ then it satisfies
conditions of Theorem 2.4 in \cite{AG0212306}. This implies $c \ge a + d +
\varepsilon$, 
where $\varepsilon = 0$ for (2),  $\varepsilon =
1$ for (3), $\varepsilon = 2$ for (4). Then, as $c \theta^2 + (d-a)\theta -
b = 0$,
$$
|\theta - \theta'|^2 = \frac{(d-a)^2 + 4 bc}{c^2} = \frac{(d+a)^2 - 4}{c^2} \le
\frac{(d+a)^2 - 4}{(d+a)^2} < 1. 
$$

Let us prove that (1) implies (2),(3) and (4). Namely, we are going to
show that (i) implies that for every $\varepsilon \le 2$ there exist $g \in
S_{\theta}$ such
that  $c > a + d + \varepsilon$. This will imply (2) for $\varepsilon = 1$, (3)
and (4) for $\varepsilon =
2$ due to Theorem 2.4 in \cite{AG0212306}.

Take any $g \in S_{\theta}$. Now, as $g$ stabilizes $\theta$, we have norm and
trace
$$
N(\theta)=-\frac b c = \frac{1 - a d}{c^2}, \; Tr(\theta) = \frac{a-d}c,
$$
and
$$
(a+d)^2 = (a-d)^2+4 a d = c^2(Tr(\theta)^2 - 4 N(\theta))+4 =
c^2|\theta-\theta'|^2+4.
$$
So, as $|\theta-\theta'| < 1$ we have $(a+d)^2 < (c-\varepsilon)^2$ if $c$ is
large enough, and $a+d < c - \varepsilon$, because $a+d>2$ and $\varepsilon \le
2$ and $c>0$. Then one can take $g^n$, which also belongs to $S_{\theta}$,
instead
of $g$, and get large enough number $c$ in the last identity.
\end{proof}

Now we are going to construct another ring, which also uses inner products in
imprimitivity bimodules $E(g^n,\theta)$. We will use left
$A_{\theta}$-valued inner products, but the same construction can be done for
the right ones. We put $R_n = \Img \lprod {A_{\theta}} {\cdot }
{\cdot} \Big|_{E_n}$~--- the vector space of finite sums of values of left inner
product on pairs of vectors from $E_n \subset E(g^n,\theta)$. In
Introduction we defined for $\Omega \in \Siegel 2$ and function $f : \Z^2
\rightarrow \Co$ periodic w.r.t. some cofinite lattice in $\Z^2$ an element
$$
\Theta[f](\Omega) = \sum_{\vec{m} \in \Z^2} f(\vec m)e(\frac 12 \vec m ^t
\Omega \vec m)e(-\frac{\theta}2 m_1 m_2)U^{m_1}V^{m_2} \in A_{\theta}.
$$

\begin{proposition}\label{quantum_theta}
$R_n = \left\{ \Theta [f](\frac 1{c_n(c_n\theta + d_n)} \Omega) \Big|
f:\Z^2\slash c_n\Z^2 \rightarrow \Co \right\}$
where 
$$
\Omega = \frac{i}{2 \Im \tau}\begin{pmatrix}|\tau|^2&-\Re
\tau\\-\Re \tau & 1 \\\end{pmatrix} \in \Siegel 2.
$$
\end{proposition}
\begin{proof}
By routine computation we get
$$
\lprod{A_{\theta}}{\phi_f}{\phi_g} = \frac 1{2(\Im \mu_n)}\sum_{\vec m \in \Z^2}
Q(\vec m)e\left(\frac12 \vec m^t \frac{\Omega}{c_n(c_n \theta
+ d_n)} \vec m\right)e(-\frac{\theta}2m_1m_2)
U^{m_1}V^{m_2}
$$
where
$$
Q(\vec m)=\sum_{\alpha \in \Z \slash c_n \Z}f(\alpha + a_n
m_1)\bar{g}(\alpha)e(\frac {\alpha}{c_n}m_2).
$$
Now the statement follows. We have $\frac 1{c_n(c_n\theta + d_n)} \Omega \in
\Siegel 2$ since $\Omega \in \Siegel 2$ and $c_n(c_n \theta + d_n) > 0$ as $g
\in
S_{\theta}$. 
\end{proof}

Note, that $R_n$ is a vector space. $\dim R_n = c_n^2 = (\dim E_n)^2$, what
implies in particular that there are no linear relations among $\lprod
{A_{\theta}} {\phi_{f_i}}{\phi_{f_j}}$ for any basis $\{ f_i \}$ in
space of functions on $\Z \slash c_n \Z$.

Now we define an operation $\underset{n,m}\star: R_n \otimes R_m \rightarrow
R_{n+m}$:
$$
\sum_i \; \lprod{A_{\theta}}{x_i}{y_i} \; \underset{n,m}\star \; \sum_j \;
\lprod{A_{\theta}}{z_j}{t_j} := \sum_{i,j} \; \lprod{A_{\theta}}{
x_i * z_j}{y_i*t_j}
$$
This operation is well defined. Indeed, every element of $R_n$ can be uniquely
represented as a linear combination of $\lprod {A_{\theta}}
{\phi_{f_i}}{\phi_{f_j}}$ as we remarked above. We can now introduce the ring $R
= \oplus_{n \ge 1} R_n$ with multiplication given by $\phi * \psi := \phi
\underset{n,m}\star \psi \in R_{n+m}$ for $\phi \in R_n$, $\psi \in R_m$.
Multiplication is obviously associative, because it is associative in the ring
$E$ defined above.
Analogously to Theorem~\ref{good_E} we have:

\begin{theorem}\label{fin_gen}
Let $\theta \in \R \backslash \Q$ be a quadratic irrationality, $\theta^{'}$ be
its Galois conjugate and $|\theta - \theta^{'}|<1$. Then there exist such $g
\in S_{\theta}$ such that the graded ring $R = R(g,\theta)$ is generated by
$R_1$ over $\Co$.
\end{theorem}
\begin{proof}
By Theorem~\ref{good_E} we can find $g \in S_{\theta}$ such that $E=E(g,\theta)$
is generated by $E_1$. So,
if we choose some basis $x_1,\dots,x_c$ in $E_1$, then $E_n$ is spanned by the
elements $x_{ i_1} * \dots * x_{i_n}$. Thus $R_n$ is spanned by elements
$$
\lprod {A_{\theta}} {x_{ i_1} * \dots * x_{i_n}}{x_{ j_1} * \dots * x_{j_n}} =
\Pi_s \; \lprod {A_{\theta}} {x_{ i_s}}{x_{j_s}}
$$
where $\lprod {A_{\theta}} {x_{ i_s}}{x_{j_s}} \in R_1$.
\end{proof}

\end{document}